\UseRawInputEncoding
\documentclass[11pt, reqno]{amsart}
\usepackage{enumitem}
\usepackage{amsthm,amsmath,amssymb}

\usepackage{mathrsfs}

\numberwithin{equation}{section}

\newtheorem{theorem}{Theorem}[section]

\theoremstyle{plain}

\newtheorem{lemma}[theorem]{Lemma}
\newtheorem{corollary}[theorem]{Corollary}
\newtheorem{remark}[theorem]{Remark}

\def\be{\begin{equation}}
	\def\ee{\end{equation}}

\usepackage{hyperref}
\hypersetup{hidelinks,
	colorlinks=true,
	allcolors=black,
	pdfstartview=Fit,
	breaklinks=true}

\numberwithin{equation}{section}
\begin{document}
\title[]{Li-Yau gradient estimates on closed manifolds under Bakry-\'Emery Ricci curvature  conditions}
\author{XingYu Song, Ling Wu }
\address{School of Mathematical Sciences and Shanghai Key Laboratory of PMMP, East China Normal University, Shanghai 200241, China}
\email{ 52215500013@stu.ecnu.edu.cn, 52215500012@stu.ecnu.edu.cn}
\date{}

\begin{abstract}
	In this paper, motivated by the work of Qi S. Zhang in \cite{ZQ}, we derive Li-Yau gradient bounds for positive solutions of the $f$-heat equation  on closed  manifolds with Bakry-\'Emery Ricci curvature bounded  below.

\end{abstract}
\maketitle

\section{Introduction}
Let $(M^n,g,e^{-f}dv)$ be a complete smooth metric measure space, where $(M^n,g)$ is an $n$-dimensional complete Riemannian manifold, $dv$ is the volume element of $g$, $f$ is a smooth function on $M$ (called the potential function), and $e^{-f}dv$ (for short, $d\mu$) is called the weighted volume element. The 
$m$-Bakry-\'Emery Ricci curvature (\cite{BE}, \cite{Lo}) associated to $(M^n,g,e^{-f}dv)$ is defined by 
$$Ric^{m,n}_f\ :=Ric+Hessf-\frac{1}{m-n}df\otimes df \ (m>n),$$
where $Ric$ is the Ricci curvature of $(M^n,g)$, $Hess$ is the Hessian with respect to the metric $g$. The $m$-Bakry-\'Emery Ricci curvature is a natural generalization of Ricci curvature on Riemannian manifolds. When $m=\infty$, we denote 
$$Ric_f=Ric^{\infty,n}_f=Ric+Hessf,$$
which is called the Bakry-\'Emery Ricci curvature (\cite{BE}). Manifolds with constant Bakry-\'Emery Ricci curvature are so called gradient Ricci solitons, which play a crucial role in the singularity analysis of the Ricci flow (\cite{Bam}, \cite{CW}, \cite{Per}, \cite{TZ}). With respect to the weighted volume element the natural self-adjoint Laplacian operator is the $f$-Laplacian 
\begin{equation}
	\Delta_f=\Delta-\left<\nabla f, \nabla \right> .\nonumber
\end{equation}
The $f$-heat equation is defined as
\begin{equation}
	( \Delta_f-\partial_t)u=0.\nonumber
\end{equation}

In \cite{LY}, Li and Yau showed that if $(M^n,g)$ is a complete Riemannian manifold with $Ric \ \ge -K$ for some constant $K \ge 0$, then for any positive solution $u$ of the heat equation $( \Delta-\partial_t)u=0$, we have 
\be\label{equ1.2}
\frac{|\nabla u|^2}{u^2}- \alpha \frac{\partial_t u}{u}\le \frac{n\alpha ^2}{2t}+\frac{n\alpha ^2K}{2(\alpha-1)}, \ \forall \alpha>1, \ t>0.
\ee
In particular, when $Ric \ \ge 0$, \ one obtains the optimal Li-Yau  bound
\be\label{equ1.3}
\frac{|\nabla u|^2}{u^2}- \frac{\partial_t u}{u}\le \frac{n}{2t}.
\ee

Many applications of \eqref{equ1.2} and \eqref{equ1.3} have been demonstrated, including the parabolic Harnack inequality, optimal Guassian estimates of the heat kernel, estimates of eigenvalues of the Laplace operator, and estimates of the Green's function. Moreover, \eqref{equ1.2} and \eqref{equ1.3} can even imply the Laplacian Comparison Theorem (see e.g. \cite{BLN} page 394).

The estimate \eqref{equ1.3} is sharp since the equality is achieved by the heat kernel of $\mathbb{R}^n$. However, \eqref{equ1.2} is not sharp for $K>0$. An open question asks if we can find sharp Li-Yau-type gradient estimates for $K>0$. Many works were done to improve or generalize \eqref{equ1.2}.

In \cite{HA}, Hamilton discovered Li-Yau-type bound for the heat equation
\begin{equation}\label{a}
	\frac{|\nabla u|^2}{u^2}-e^{2Kt}\frac{\partial_t u}{u}\le e^{4Kt}\frac{n}{2t}.
\end{equation} 

In \cite{Yau}, Yau obtained the following estimate 
\begin{equation}\label{b}
	\frac{|\nabla u|^2}{u^2}-\frac{\partial_t u}{u}\le \frac{n}{2t}+\sqrt{2nK}\sqrt{\frac{|\nabla u|^2}{u^2}+\frac{n}{2t}+2nK} .
\end{equation}

In \cite{BQ}, Bakry and Qian obtained
\begin{equation}\label{c}
	\frac{|\nabla u|^2}{u^2}-\left(1+\frac{2}{3}Kt\right)\frac{\partial_t u}{u}\le \frac{n}{2t}+\frac{nK}{2}\left(1+\frac{K}{3}t\right).
\end{equation}

In \cite{Qian Bin}, Bin Qian improved the estimate \eqref{c} in the following form
\begin{equation}\label{d}
	\frac{|\nabla u|^2}{u^2}-\left(1+\frac{2K}{a(t)}\int_{0}^{t}a(s)ds\right)\frac{\partial_t u}{u}\le \frac{nK}{2}+\frac{nK^2}{2a(t)}\int_{0}^{t}a(s)ds+\frac{n}{8a(t)}\int_{0}^{t}\frac{a'(s)^2}{a(s)}ds,
\end{equation}
where $a(t):(0,\infty)\rightarrow (0,\infty)$ is any $C^1$ positive function that satisfies the following two conditions:

(A1) For all $t>0$, $a(t)>0$, $a'(t)>0$ and $\lim\limits_{t \rightarrow 0}a(t)=\lim\limits_{t \rightarrow 0}\frac{a(t)}{a'(t)}=0$.

(A2) For any $L>0$, $\frac{a'(t)^2}{a(t)}$ is continuous and integrable on the interval $[0,L]$.

In the above results, the constant $\alpha$ in \eqref{equ1.2} is replaced by functions in the form of $\alpha(t,K)$ which is strictly greater than $1$ but converges to $1$ as $t\rightarrow 0$. Recently, Qi S. Zhang \cite{ZQ} obtains that for all closed manifolds one can take $\alpha=1$ for $K\ \ge 0$.
It shows that if $(M^n,g)$ is a closed Riemannian manifold with $Ric \ \ge -K$ for $K \ge 0$, and $diam_{M}$ is  the diameter of $M$,\ then for any positive solution $u$ of the heat equation $( \Delta-\partial_t)u=0$,
\begin{equation}\label{equ1.4}
	t(\frac{|\nabla u|^2}{u^2}-\frac{\partial_t u}{u})\le \frac{n}{2}+\sqrt{2nK(1+Kt)(1+t)}diam_M+ \sqrt{K(1+Kt)(C_1t+C_2Kt)},
\end{equation}
where  $C_1$ and $C_2$ are positive constants only depending on $n$.
It is an improvement on  Li-Yau gradient bound (\ref{equ1.2}) on closed manifolds. 
 
In (\cite{ZZ}, \cite{ZZ2}), Qi S. Zhang and Meng Zhu obtained Li-Yau-type gradient estimates under integral curvature assumptions. Moreover, Li-Yau-type bounds were also got for weighted manifolds with Bakry-\'Emery Ricci curvature bounded below \cite{Li}. More information about Li-Yau-type  bounds can be found in (\cite{GM}, \cite{LX}, \cite{QZZ}, \cite{WFY}, \cite{ZZH}, \cite{ZFF}, \cite{ZFF2}, \cite{ZFF3}).

In this paper, we show the optimal Li-Yau gradient bounds for $f$-heat equation on closed manifolds with either  Bakry-\'Emery Ricci curvature or $m$-Bakry-\'Emery Ricci curvature bounded below, which both generalize \eqref{equ1.4}.

More precisely, we show that
 \begin{theorem}\label{thm1} 
 Let $(M^n,g)$ be a  closed Riemannian manifold with  $Ric_f \ge -K$ and $|\nabla f| \le L$ for some constants $K, L \ge 0$. Let $u$ be a positive solution of the f-heat equation on $M \times (0,+\infty)$, i.e., $( \Delta_f-\partial_t)u=0$, and $diam_{M}$  the diameter of $M$. Then there exists a contant $c$ depending only on $n$, such that 
 \begin{equation}\label{1.5}
 	\begin{aligned}
 		t(\frac{|\nabla u|^2}{u^2}-\frac{\partial_t u}{u})\le &\frac{n}{2}+c(L+\sqrt{K})\sqrt{(1+Kt)(1+t)}diam_M\\
 		&+ c(L+\sqrt{K})\sqrt{(1+Kt)(t+Kt+L^2t+A^2t+A^2Kt)},
 	\end{aligned} 
 \end{equation}
 where $A=\sup\limits_{x\in M}|f(x)|$.

\end{theorem}
Notice that (\ref{1.5}) still holds when we add $f$ by any constant. In particular, for some point $o\in M$ we may choose $f$ such that $f(o)$ equals zero. Then (\ref{1.5}) becomes 
\begin{equation}\label{1.6}
	\begin{aligned}
		t(\frac{|\nabla u|^2}{u^2}-\frac{\partial_t u}{u})\le &\frac{n}{2}+c(L+\sqrt{K})\sqrt{(1+Kt)(1+t)}diam_M\\
		&+ c(L+\sqrt{K})\sqrt{(1+Kt)(t+Kt+L^2t+L^2diam_M^2t+L^2diam_M^2Kt)}.
	\end{aligned}
\end{equation}

If the potential function is constant, then from (\ref{1.6}) we get a result similar to (\ref{equ1.4}).

\begin{corollary}\label{cor1}
	Let $(M^n,g)$ be a  closed Rimannian manifold with $Ric \ge -K$ where $K \ge 0$.  Let $u$ be a positive solution of the heat equation on $M \times (0,+\infty)$, i.e., $( \Delta-\partial_t)u=0$, and $diam_{M}$  the diameter of $M$.  Then there exists
	a constant $c$ depending only on $n$, such that 
	 \begin{equation}\label{1.7}
	 	t(\frac{|\nabla u|^2}{u^2}-\frac{\partial_t u}{u})\le \frac{n}{2}+c\sqrt{K(1+Kt)(1+t)}diam_M+  c\sqrt{K(1+Kt)(t+Kt)}.
	 \end{equation}
\end{corollary}

For closed manifolds with $m$-Bakry-\'Emery Ricci curvature bounded below, we have a similar result. But there is no assumption on the potential function and the constants depend on $m$.

\begin{theorem}\label{thm2}
	 Let $(M^n,g)$ be a  closed Riemannian manifold with  $Ric^{m,n}_f \ge -K$\ for some constant $K \ge 0$. Let $u$ be a positive solution of the f-heat equation on $M \times (0,+\infty)$, i.e., $( \Delta_f-\partial_t)u=0$, and $diam_{M}$  the diameter of $M$. Then there exist contants $\widetilde{c_1} \ and \ \widetilde{c_2}$ depending only on $m$ such that 
	  \begin{equation}\label{1.8}
	  	t(\frac{|\nabla u|^2}{u^2}-\frac{\partial_t u}{u})\le \frac{m}{2}+\sqrt{2mK(1+Kt)(\widetilde{c_1}+t)}diam_M+  \widetilde{c_2}\sqrt{K(1+Kt)(t+Kt)}.
	  \end{equation}

\end{theorem}

If the potential function is constant, then we can take $m=n$ and get a result similar to  Corollary \ref{cor1}.
\begin{remark}
	The constants  $c$, $\widetilde{c_1} \ and \ \widetilde{c_2}$ in Theorem \ref{thm1} and Theorem \ref{thm2} arise from the volume comparison theorem and  upper and lower bounds for the $f$-heat kernel. Since the lower bound constants of the $f$-heat kernel cannot be written accurately, the constants $c$\and \ $\widetilde{c_1}$ cannot be written as  accurately as Qi S. Zhang does.
\end{remark}

We prove Theorem \ref{thm1} and Theorem \ref{thm2} separately in sections 2 and 3. The proofs follow a method of Qi S. Zhang \cite{ZQ}. It is sufficient to prove the optimal Li-Yau bound for the $f$-heat kernel, then the same bound holds for all positive solutions of the $f$-heat equation. The main tools in the proof of Theorem \ref{thm1} are relative volume comparison \cite{WWW} of Guofang Wei and  Will Wylie, Yi Li's Hamilton type estimates \cite{Li}, the upper and lower bounds for the $f$-heat kernel \cite{WWu} and the Harnack inequality for the positive solutions of the $f$-heat equation \cite{WWu} proved by Jiayong Wu and Peng Wu for the manifolds with Bakry-\'Emery Ricci curvature bounded below. Here we generally assume that the potential function is bounded or the gradient of the function is bounded. The proof of Theorem \ref{thm2} is similar to Theorem \ref{thm1}. However, in Theorem \ref{thm2} for manifolds with $m$-Bakry-\'Emery Ricci curvature bounded below, there is no assumption on potential function and relative results can be found in (\cite{ALR}, \cite{Li}, \cite{WWW}).

\section{ Li-Yau gradient bounds on closed manifolds under  Bakry-\'Emery Ricci curvature conditions  }
In this section, we prove Theorem \ref{thm1}. First we present some results as preparation.
\begin{lemma} \label{lem1} (\cite{Li})(Hamilton type estimate)
	Suppose that $(M^n,g)$ is a closed Riemannian manifold with $Ric_f \ge -K $ where $K\ge 0$. If $u$ is a positive solution of $( \Delta_f-\partial_t)u=0$ with $0<u\le B$ on $M \times (0,T]$ for some constant $B$, then 
	\begin{equation}
		\frac{|\nabla u|^2}{u^2} \le (\frac{1}{t}+2K)\ln\frac{B}{u},\nonumber
	\end{equation}
	on $M \times (0,T]$.
\end{lemma}

In \cite{WWW}, Guofang Wei and Will Wylie proved the volume comparison theorem under Bakry-\'Emery Ricci curvature conditions.
\begin{lemma}\label{lem2} (\cite{WWW}) (volume comparison)
Let $(M^n,g,e^{-f}dv)$ be a complete smooth metric measure space with $Ric_f \ge(n-1)H$. For $p\in M$,
If $|f(x)| \le A$ for some constant $A$, then for $R\ge r>0$ $(assume \  R \le \frac{\pi}{4\sqrt{H}}\ if \ H>0)$, 
\[\frac{V_f(B_p(R))}{V_f(B_p(r))} \le \frac{V^{n+4A}_H(B(R))}{V^{n+4A}_H(B(r))},\]
where $V^{n}_H(B(r))$ is the volume of the geodesic ball with radius $r$ in the model space $M^n_H$ and $V_f(B_p(r))$ is the weighted volume of the geodesic ball with radius $r$ in $M$ centerd at $p$.
\end{lemma}

In \cite{WWu}, Jiayong Wu and Peng Wu showed that
\begin{lemma}\label{lem3} (\cite{WWu})
Let $(M^n,g,e^{-f}dv)$ be a complete noncompact smooth metric measure space with $Ric_f \ge -(n-1)K$ for $K \ge 0$. For any point $o\in M$ and $R>0$, denote $A(R)=\sup\limits_{x\in B_o(3R)}|f(x)|,\ A^{'}(R)=\sup\limits_{x\in B_o(3R)}|\nabla f(x)|$, and let $H(x,t,y)$ be the f-heat kernel. 
Then for any $\epsilon > 0$, there exist constants $c_3(n,\epsilon),c_i(n)\ 4\le i \le 8$ such that 
\begin{equation}\label{2.1}
	\frac{c_3e^{c_4A+c_5(1+A)\sqrt{Kt}}}{V_f^{\frac{1}{2}}(B_x(\sqrt{t}))V_f^{\frac{1}{2}}(B_y(\sqrt{t}))}e^{\left(-\frac{d^2(x,y)}{(4+\epsilon)t}\right)} \ge H(x,t,y) \ge \frac{c_6e^{-c_7(A'^2+K)t}}{V_f(B_x(\sqrt{t}))}e^{\left(-\frac{d^2(x,y)}{c_8t}\right)}
\end{equation}
for all $x,\ y \in B_o(\frac{1}{2}R)$ and $0<t<\frac{R^2}{4}$, where $d(x,y)$ is the geodesic distance between $x$ and $y$.
	
\end{lemma}
Note that for closed Riemannian manifolds with $Ric_f \ge -K$ and  $|\nabla f| \le L$ where $K,\ L\ge0$, the same bound still holds for all $x,\ y\in M$ and $t>0$.

\begin{lemma}\label{lem4} (\cite{WWu})
	Let $(M^n,g,e^{-f}dv)$ be a complete noncompact smooth metric measure space with $Ric_f \ge -(n-1)K$ for $K \ge 0$. For any point $o\in M$ and $R>0$, denote $A^{'}(R)=\sup\limits_{y\in B_o(R+1)}|\nabla f(y)|$. There exists a constant $c(n)$ such that, for two positive solutions $u(x,s)$ and $u(y,t)$ of the f-heat equation in $B_o\left(\frac{R}{2}\right) \times (0,T),\ 0<s<t<T$,
\[\ln \left(\frac{u(x,s)}{u(y,t)}\right) \le c(n)\left[\left(A'^2+K+\frac{1}{R^2}+\frac{1}{s}\right)(t-s)+\frac{d^2(x,y)}{t-s}\right].\]
	
\end{lemma}
Note that for closed Riemannian manifolds with $Ric_f \ge -K$ and  $|\nabla f| \le L$ where $K,\ L\ge0$, there exists a constant $c(n)$ such that, in $M \times (0,+\infty),\ 0<s<t$, 
\begin{equation}\label{2.2}
	u(x,s)\le u(y,t)e^{c(n)\left[\left(L^2+K+\frac{1}{s}\right)(t-s)+\frac{d^2(x,y)}{t-s}\right]}.
\end{equation}

Now we are ready to prove Theorem \ref{thm1}.

\noindent{\it Proof of Theorem \ref{thm1}.}
Let $u=u(x,t)$ be a positive solution of the $f$-heat equation on $M\times(0,+\infty)$. Denote
\begin{equation}
	Y=Y(x,t)=|\nabla \ln u|^2-\partial_t(\ln u)=\left(\frac{|\nabla u|^2}{u^2}-\frac{\partial_t u}{u}\right)=-\Delta_f \ln u .\nonumber
\end{equation}
The Bochner formula \cite{Li} for  $Ric_f$ is 
\begin{equation} \label{equ2.3}
	\begin{aligned}
		\frac{1}{2}\Delta_f|\nabla u|^2=|Hess\ u|^2+\left<\nabla \Delta_f u,\nabla u\right>+Ric_f(\nabla u,\nabla u).
	\end{aligned}
\end{equation}
Using \eqref{equ2.3} and the Bakry-\'Emery Ricci curvature condition, we have 
\begin{equation}\label{eq1}
\begin{aligned}
	(\Delta_f-\partial_t)Y+&2\left<\nabla Y,\nabla \ln u\right>=2|Hess \ln u|^2+2Ric_f(\nabla \ln u,\nabla \ln u)\\
	&\ge \frac{2(\Delta \ln u )^2}{n}-2K|\nabla \ln u |^2\\
	&=\frac{2(\Delta_f \ln u+\left<\nabla f,\nabla \ln u\right> )^2}{n}-2K|\nabla \ln u |^2\\
	&=\frac{2}{n}(\Delta_f \ln u)^2+\frac{4}{n}(\Delta_f \ln u)\left<\nabla f,\nabla \ln u\right>+\frac{2}{n}|\left<\nabla f,\nabla \ln u\right>|^2-2K|\nabla \ln u|^2\\
	&\ge \frac{2}{n}Y^2-\frac{4}{n}\left<\nabla f,\nabla \ln u\right>Y-\left(\frac{2}{n}L^2+2K\right)|\nabla \ln u|^2.
\end{aligned}
\end{equation}
Let $Y^+(x,t)=\max\{Y(x,t),0\}$. Then inequality \eqref{eq1} implies that $Y^+$ is a subsolution of the inequality in the weak sense: on $M \times(0,+\infty)$,
\begin{equation} \label{eq2}
	(\Delta_f-\partial_t)Y^+ +2\left<\nabla Y^+,\nabla \ln u \right> \ge \frac{2}{n}(Y^+)^2-\frac{4}{n}\left<\nabla f,\nabla \ln u\right>Y^+-\left(\frac{2}{n}L^2+2K\right)|\nabla \ln u|^2.
\end{equation}
For a positive integer $j$ and a small positive number $\epsilon$, multiplying both sides of \eqref{eq2} by $((t-\epsilon)^+)^{2j+2}(Y^+)^{2j}$ and integrating on $M\times (0,T]$ for any $T>0$ give that
\begin{equation}\label{eq3}
	\begin{aligned}
		&\frac{2}{n}\int_{0}^{T}\int_M ((t-\epsilon)^+Y^+)^{2j+2} d\mu dt \le \int_{0}^{T}\int_M ((t-\epsilon)^+)^{2j+2}(Y^+)^{2j}(\Delta_f -\partial_t)Y^+ d\mu dt\\
		&+ 2\int_{0}^{T}\int_M ((t-\epsilon)^+)^{2j+2}(Y^+)^{2j} \left<\nabla Y^+,\nabla \ln u \right>d\mu dt\\
		&+ (\frac{2}{n}L^2+2K)\int_{0}^{T}\int_M |\nabla \ln u|^2((t-\epsilon)^+)^{2j+2}(Y^+)^{2j}d\mu dt\\
		&+\frac{4L}{n}\int_{0}^{T}\int_M |\nabla \ln u|((t-\epsilon)^+)^{2j+2}(Y^+)^{2j+1}d\mu dt=T_1+T_2+T_3+T_4.
	\end{aligned} 
\end{equation}
Here we need to estimate the upper bound of $(t-\epsilon)^{+}Y^+$ on $M \times (0,T]$. Then let us bound $T_1$, $T_2$, $T_3$ and $T_4$ respectively. Using integration by parts, we see that
\begin{equation}\label{eq4}
	\begin{aligned}
		T_1&=-2j\int_{0}^{T}\int_M ((t-\epsilon)^+)^{2j+2}(Y^+)^{2j-1}|\nabla Y^+|^2d\mu dt -\frac{1}{2j+1}\int_M \int_{0}^{T} ((t-\epsilon)^+)^{2j+2}d(Y^+)^{2j+1} d\mu \\
		&=-2j\int_{0}^{T}\int_M ((t-\epsilon)^+)^{2j+2}(Y^+)^{2j-1}|\nabla Y^+|^2d\mu dt-\frac{1}{2j+1}\int_M ((t-\epsilon)^+)^{2j+2}(Y^+)^{2j+1}|^{T}_0 d\mu \\
		&+\frac{2j+2}{2j+1}\int_{0}^{T}\int_M ((t-\epsilon)^+Y^+)^{2j+1} d\mu dt \\
		&=-2j\int_{0}^{T}\int_M ((t-\epsilon)^+)^{2j+2}(Y^+)^{2j-1}|\nabla Y^+|^2d\mu 
		dt-\frac{1}{2j+1}\int_M ((T-\epsilon)^+)^{2j+2}(Y^+)^{2j+1}(x,T) d\mu  \\
		&+\frac{2j+2}{2j+1}\int_{0}^{T}\int_M ((t-\epsilon)^+Y^+)^{2j+1} d\mu dt.
	\end{aligned}
\end{equation}
Writing $(Y^+)^{2j} \nabla Y^+ =\frac{1}{2j+1}\nabla(Y^+)^{2j+1}$ and doing integration by parts, we deduce
\begin{equation} \label{eq5}
	\begin{aligned}
		T_2&=\frac{2}{2j+1} \int_{0}^{T}\int_M \left<\nabla (Y^+)^{2j+1},\nabla \ln u \right> ((t-\epsilon)^+)^{2j+2} d\mu dt\\
		&=-\frac{2}{2j+1} \int_{0}^{T}\int_M (\Delta_f \ln u)(Y^+)^{2j+1} ((t-\epsilon)^+)^{2j+2}d\mu dt \\
		&\le \frac{2}{2j+1} \int_{0}^{T}\int_M ((t-\epsilon)^+Y^+)^{2j+2} d\mu dt.
	\end{aligned}
\end{equation}
Throwing away the non-positive terms of \eqref{eq4} and plugging (\ref{eq4}) and (\ref{eq5}) into (\ref{eq3}) give that 
\begin{equation} \label{eq6}
	\begin{aligned}
		&(\frac{2}{n}-\frac{2}{2j+1})\int_{0}^{T}\int_M ((t-\epsilon)^+Y^+)^{2j+2}d\mu dt\\
		&\le \frac{2j+2}{2j+1}\int_{0}^{T}\int_M ((t-\epsilon)^+Y^+)^{2j+1} d\mu dt +\underbrace{(\frac{2}{n}L^2+2K)\int_{0}^{T}\int_M |\nabla \ln u|^2((t-\epsilon)^+)^{2j+2}(Y^+)^{2j}d\mu dt}_{T_3}\\
		&+\underbrace{\frac{4}{n}L\int_{0}^{T}\int_M |\nabla \ln u|((t-\epsilon)^+)^{2j+2}(Y^+)^{2j+1}d\mu dt}_{T_4}.
	\end{aligned}
\end{equation}

This estimate holds for all positive solutions. Now we take, in particular $u=H(x,t,y)$ the $f$-heat kernel with pole at a fixed point $y\in M$. We will find upper bounds for $T_3$ and $T_4$, which rely on the Hamilton type estimate (Lemma \ref{lem1}).

For a time $t_0>0$, we consider the $f$-heat kernel $H(x,t+t_0,y)$ with $t\in[0,t_0]$.

According to Lemma \ref{lem3} , we can choose $\epsilon=1$ and find the upper and lower bounds for $f$-heat kernel. There exist some positive constants $C_1, C_2, C_3, C_4, C_5$ and $ C_6$ depending only on $n$, such that
\begin{equation}\label{eq7}
\frac{C_1e^{C_2A+C_3(1+A)\sqrt{K(t+t_0)}}}{V_f^{\frac{1}{2}}(B_x(\sqrt{t+t_0}))V_f^{\frac{1}{2}}(B_y(\sqrt{t+t_0}))}e^{\left(-\frac{d^2(x,y)}{5(t+t_0)}\right)} \ge H(x,t+t_0,y) \ge \frac{C_4e^{-C_5(L^2+K)(t+t_0)}}{V_f(B_x(\sqrt{t+t_0}))}e^{\left(-\frac{d^2(x,y)}{C_6(t+t_0)}\right)}.
\end{equation}
This upper bound implies
\begin{equation}
	B:=\sup_{ M\times(0,t_0)}H(x,t+t_0,y)\le \frac{C_1e^{C_2A+C_3(1+A)\sqrt{2Kt_0}}}{\inf_{x\in M}V_f(B_x(\sqrt{t_0}))},\nonumber
\end{equation}
which yields, by the lower bound of  $H(x,t+t_0,y)$, that
\begin{equation}\label{eq8}
	\frac{B}{H(x,t+t_0,y)}\le C_7e^{\left(C_2A+C_3(1+A)\sqrt{2Kt_0}+C_5(L^2+K)(2t_0)+ \frac{d^2(x,y)}{C_6t_0}\right)} 
	\frac{\sup_{y\in M}V_f(B_y(\sqrt{2t_0}))}{\inf_{x\in M}V_f(B_x(\sqrt{t_0}))}.
\end{equation} 
We notice that the infimum and supremum of the geodesic ball's volumes can be achieved by some points on $M$, say point $p$ and point $q$, i.e.,
\begin{equation}\label{eq9}
\frac{\sup_{y\in M}V_f(B_y(\sqrt{2t_0}))}{\inf_{x\in M}V_f(B_x(\sqrt{t_0}))}=\frac{V_f(B_q(\sqrt{2t_0}))}{V_f(B_p(\sqrt{t_0}))}.
\end{equation}
Here is the ratio of geodesic balls's volumes at different points. Notice that the relationship between the geodesic ball $B_q(\sqrt{2t_0})$ and the geodesic ball $B_p(\sqrt{2t_0}+d(p,q))$ is $B_q(\sqrt{2t_0})\subset B_p(\sqrt{2t_0}+d(p,q))$. Then  applying the volume comparison theorem (Lemma \ref{lem2}), we get
\begin{equation}\label{eq10}
	\begin{aligned}
		\frac{V_f(B_p(\sqrt{2t_0}+d(p,q)))}{V_f(B_p(\sqrt{2t_0}))}
		&\le \frac {\int_{0}^{d(p,q)+\sqrt{2t_0}}\left(sinh\left(\sqrt{\frac{K}{n-1}}r\right)\right)^{n+4A-1}dr}{\int_{0}^{\sqrt{2t_0}}\left(sinh\left(\sqrt{\frac{K}{n-1}}r\right)\right)^{n+4A-1}dr}\\ &\le \left(\frac{d(p,q)}{\sqrt{2t_0}}+1\right)^{n+4A}e^{d(p,q) \sqrt{\frac{K}{n-1}}(n+4A-1)}\\
		&=e^{(n+4A) \ln (\frac{d(p,q)}{\sqrt{2t_0}}+1) +d(p,q)\sqrt{\frac{K}{n-1}}(n+4A-1)}\\
		&\le e^{(n+4A)  \frac{d(p,q)}{\sqrt{2t_0}}+ d(p,q)\sqrt{\frac{K}{n-1}}(n+4A-1)},  
	\end{aligned}
\end{equation}
and

\begin{equation}\label{eq11}
\begin{aligned}
	\frac{V_f(B_p(\sqrt{2t_0}))}{V_f(B_p(\sqrt{t_0}))} 
	&\le \frac {\int_{0}^{\sqrt{2t_0}}\left(sinh\left(\sqrt{\frac{K}{n-1}}r\right)\right)^{n+4A-1}dr}{\int_{0}^{\sqrt{t_0}}\left(sinh\left(\sqrt{\frac{K}{n-1}}r\right)\right)^{n+4A-1}dr}\\ &\le
	2^{\left(\frac{n+4A}{2}\right)}e^{[(\sqrt{2}-1)\sqrt{t_0}\sqrt{\frac{K}{n-1}}](n+4A-1)}.  
\end{aligned}
\end{equation}
Substituting (\ref{eq10}) and (\ref{eq11}) into (\ref{eq9}), we find
\begin{equation}
\begin{aligned}
	\frac{\sup_{y\in M}V_f(B_y(\sqrt{2t_0}))}{\inf_{x\in M}V_f(B_x(\sqrt{t_0}))}&=\frac{V_f(B_q(\sqrt{2t_0}))}{V_f(B_p(\sqrt{t_0}))}  \le\frac{V_f(B_p(\sqrt{2t_0}+d(p,q)))}{V_f(B_p(\sqrt{t_0}))}\nonumber\\
	&=\frac{V_f(B_p(\sqrt{2t_0}+d(p,q)))}{V_f(B_p(\sqrt{2t_0}))} \frac{V_f(B_p(\sqrt{2t_0}))}{V_f(B_p(\sqrt{t_0}))}\nonumber\\
	&\le 2^{\left(\frac{n+4A}{2}\right)}e^{[(\sqrt{2}-1)\sqrt{t_0}\sqrt{\frac{K}{n-1}}](n+4A-1)+{(n+4A) \frac{d(p,q)}{\sqrt{2t_0}}+ d(p,q)\sqrt{\frac{K}{n-1}}(n+4A-1)}}\nonumber.
\end{aligned}
\end{equation}
This and (\ref{eq8}) imply that
\begin{equation}
	\begin{aligned}
		&\ln \frac{B}{H(x,t+t_0,y)} \le \ln C_7 +C_2A+C_3\sqrt{2Kt_0}+C_3A\sqrt{2Kt_0}+2C_5L^2t_0+2C_5Kt_0\\
		&+\frac{d^2(x,y)}{C_6t_0}+\frac{nd(p,q)}{\sqrt{2t_0}}+2\sqrt{2}A\frac{d(p,q)}{\sqrt{t_0}}+d(p,q)\sqrt{K(n-1)}+4d(p,q)A\sqrt{\frac{K}{n-1}}\\
		&+\sqrt{K(n-1)t_0}+4A\sqrt{\frac{Kt_0}{n-1}}+\frac{\ln 2}{2}(n+4A).
	\end{aligned} \nonumber
\end{equation}
Using the following basic inequalities with
\begin{equation}
	\sqrt{Kt_0}\le \frac{1}{4}+Kt_0,\  \frac{d(p,q)}{\sqrt{t_0}}\le \frac{d^2(p,q)}{t_0} +\frac{1}{4},\  \frac{Ad(p,q)}{\sqrt{t_0}}\le \frac{d^2(p,q)}{t_0}+\frac{1}{4}A^2,\nonumber
\end{equation}

\begin{equation}
	d(p,q)\sqrt{K}\le \frac{d^2(p,q)}{t_0}+\frac{1}{4}Kt_0,\  Ad(p,q)\sqrt{K}\le \frac{d^2(p,q)}{t_0}+\frac{1}{4}A^2Kt_0,\ \and\ A\sqrt{Kt_0}\le \frac{1}{4}A^2+Kt_0,\nonumber
\end{equation}
we get
\begin{equation}
	\ln  \frac{B}{H(x,t+t_0,y)}\le C_8\left(1+A+Kt_0+AKt_0+L^2t_0+\frac{diam_M^2}{t_0}+A^2+A^2Kt_0\right),\ \forall t\in(0,t_0].\nonumber
\end{equation} 
Using $A\le \frac{1}{4}+A^2$ gives that
\begin{equation}
	\ln  \frac{B}{H(x,t+t_0,y)}\le C_{9}\left(1+Kt_0+L^2t_0+\frac{diam_M^2}{t_0}+A^2+A^2Kt_0\right).\nonumber
\end{equation}
Hamilton type estimate (Lemma \ref{lem1}) implies that 
\begin{equation}
	t|\nabla_x \ln H(x,t+t_0,y)|^2 \le C_{9}(1+2Kt_0)\left(1+Kt_0+L^2t_0+\frac{diam_M^2}{t_0}+A^2+A^2Kt_0\right),\ t\in (0,t_0].\nonumber
\end{equation}
We take $t=t_0$ and use the arbitrariness of $t_0$ to conclude
\begin{equation} \label{eq12}
	\begin{aligned}
		t|\nabla_x \ln H(x,t,y)|^2 \le C_{10}(1+Kt)\left(1+Kt+L^2t+\frac{diam_M^2}{t}+A^2+A^2Kt\right),\ \forall t>0
	\end{aligned}.
\end{equation}

The bound is adequate for us when the time is short, say $t \le 4$. When $t$ is large, the $f$-heat kernel converges to the positive constant $\frac{1}{V_f(M)}$, where $V_f(M)$ is the weighted volume of $M$. In this case the above bound becomes inaccurate. Instead we will use a better bound based on the Harnack inequality \eqref{2.2}.

Pick any time $t\ge4$. Since $\int_M H(x,t+1,y)d\mu(x) =1$, there is a point $x_1\in M$ such that  $H(x_1,t+1,y)=\frac{1}{V_f(M)}$. According to \eqref{2.2}, there exists a dimensional constant $C_0 >0$ such that 
\begin{equation} 
	H(x,t,y) \le H(x_1,t+1,y)e^{C_0(L^2+K+\frac{1}{t}+d^2(x,x_1))}.\nonumber
\end{equation}
Since $t\ge4$, this implies
\begin{equation}\label{eq13}
	H(x,t,y) \le \frac{1}{V_f(M)}e^{C_0(L^2+K+\frac{1}{4}+diam_M^2)}:=B.
\end{equation}
Similarly, there is a point $x_2$ such that $H(x_2,t-1,y)=\frac{1}{V_f(M)}$ and that
\begin{equation} 
	H(x_2,t-1,y) \le H(x,t,y)e^{C_0(L^2+K+\frac{1}{t-1}+d^2(x,x_2))},\nonumber
\end{equation}
which infers
\begin{equation} \label{eq14}
H(x,t,y) \ge \frac{1}{V_f(M)} e^{-C_0(L^2+K+\frac{1}{3}+diam_M^2)}.
\end{equation}
Using (\ref{eq13}) and (\ref{eq14}), we find, for $t_0 \ge4$, that 
\begin{equation}
	\ln \frac{B}{H(x,t+t_0,y)} \le 2C_0(L^2+K+1+diam^2_M),\ t\in (0,t_0].\nonumber
\end{equation}
This and Lemma \ref{lem1} yield
\begin{equation}
	t|\nabla_x \ln H(x,t+t_0,y)|^2 \le 2C_0(1+2Kt_0)(L^2+K+1+diam^2_M),\ t\in (0,t_0].\nonumber
\end{equation}
Therefore
\begin{equation} \label{eq15}
	\begin{aligned}
		t|\nabla_x \ln H(x,t,y)|^2 \le 4C_0(1+Kt)(L^2+K+1+diam^2_M),\ t\ge4.
	\end{aligned}
\end{equation}
Next plugging (\ref{eq12}) for $t<4$ and (\ref{eq15}) for $t\ge4$ into the term $T_3$ and $T_4$ in (\ref{eq6}) with $u=H(x,t,y)$, we obtain
\begin{equation}
	\begin{aligned}
		T_3 &\le \left(\frac{2}{n}L^2+2K\right)\left(\int_{0}^{4}\int_M |\nabla \ln u|^2t^2((t-\epsilon)^+Y^+)^{2j}d\mu dt \right)\\
		&\ \ \ \ +\left(\frac{2}{n}L^2+2K\right)\left(\int_{4}^{T}\int_M |\nabla \ln u|^2t^2((t-\epsilon)^+Y^+)^{2j}d\mu dt\right)\\
		&\le C_{11}\left(L^2+K\right)(1+KT)(T+KT+L^2T+A^2T+A^2KT+diam^2_M\\
		&\ \ \ \ +Tdiam^2_M)\int_{0}^{T}\int_M ((t-\epsilon)^+Y^+)^{2j}d\mu dt,
	\end{aligned} \nonumber
\end{equation}
and
\begin{equation}
		\begin{aligned}
		T_4 &\le \frac{4}{n}L\int_{0}^{4}\int_M t|\nabla \ln u|((t-\epsilon)^+Y^+)^{2j+1}d\mu dt  +\frac{4}{n}L\int_{4}^{T}\int_M t|\nabla \ln u|((t-\epsilon)^+Y^+)^{2j+1}d\mu dt\\
		&\le C_{12}L\sqrt{(1+KT)(T+KT+L^2T+A^2T+A^2KT+diam^2_M+Tdiam^2_M)}\\
		&\ \ \ \times\int_{0}^{T}\int_M ((t-\epsilon)^+Y^+)^{2j+1}d\mu dt.
	\end{aligned} \nonumber
\end{equation}
Denote $\lambda=\sqrt{(1+KT)(T+KT+L^2T+A^2T+A^2KT+diam^2_M+Tdiam^2_M)}$ and (\ref{eq6}) becomes
\begin{equation} \label{eq16}
	\begin{aligned}
		( &\frac{2}{n}-\frac{2}{2j+1} )\int_{0}^{T}\int_M ((t-\epsilon)^+Y^+)^{2j+2}d\mu dt \le \left(\frac{2j+2}{2j+1}+C_{12}L\lambda\right)\int_{0}^{T}\int_M ((t-\epsilon)^+Y^+)^{2j+1}d\mu dt\\
		&+C_{11}\left(L^2+K\right)\lambda^2\int_{0}^{T}\int_M ((t-\epsilon)^+Y^+)^{2j}d\mu dt.
	\end{aligned}
\end{equation}
Using the notation
\begin{equation}
	A_{j,\epsilon}=\left(\int_{0}^{T}\int_M ((t-\epsilon)^+Y^+)^{j}d\mu dt\right) ^{\frac{1}{j}},\nonumber
\end{equation}
we can write (\ref{eq16}) as 
\begin{equation} \label{eq17}
	\begin{aligned}
		(\frac{2}{n}-\frac{2}{2j+1})A_{2j+2,\epsilon}^{2j+2} \le \left(\frac{2j+2}{2j+1}+C_{12}L\lambda\right)A_{2j+1,\epsilon}^{2j+1}+C_{11}\left(L^2+K\right)\lambda^2A_{2j,\epsilon}^{2j}.
	\end{aligned}
\end{equation}
By H\"older inequality,
\begin{equation}
	A_{2j+1,\epsilon}^{2j+1} \le A_{2j+2,\epsilon}^{2j+1}\left(\int_{0}^{T}\int_Md\mu dt\right)^{\frac{1}{2j+2}},\ A_{2j,\epsilon}^{2j} \le A_{2j+2,\epsilon}^{2j}\left(\int_{0}^{T}\int_Md\mu dt\right)^{\frac{2}{2j+2}},\nonumber
\end{equation}
which imply, together with (\ref{eq17}), that 
\begin{equation} \label{eq18}
	\begin{aligned}
		( &\frac{2}{n}-\frac{2}{2j+1} )A_{2j+2,\epsilon}^{2} \le\left(\frac{2j+2}{2j+1}+C_{12}L\lambda\right)A_{2j+2,\epsilon}\left(\int_{0}^{T}\int_Md\mu dt\right)^{\frac{1}{2j+2}}\\
		&+C_{11}\left(L^2+K\right)\lambda^2\left(\int_{0}^{T}\int_Md\mu dt\right)^{\frac{2}{2j+2}}.
	\end{aligned}
\end{equation}
Letting $j\to \infty$, we arrive at
\begin{equation}
\begin{aligned}
	\frac{2}{n}A_{\infty,\epsilon}^2\le \left(1+C_{12}L\lambda\right)A_{\infty,\epsilon}+C_{11}(L^2+K)\lambda^2,
\end{aligned}
\end{equation}
where $A_{\infty,\epsilon}=\sup\limits_{ M\times [0,T]}(t-\epsilon)^+Y^+.$

This shows 
\begin{equation}
\begin{aligned}
	\sup_{ M\times [0,T]}(t-\epsilon)^+Y^+&\le \frac{n}{2}+C_{13}(L+\sqrt{K})\sqrt{(1+KT)(1+T)}diam_M\\
	&\ \ \ + C_{13}(L+\sqrt{K})\sqrt{(1+KT)(T+KT+L^2T+A^2T+A^2KT)}.
\end{aligned}
\end{equation}
Since $\epsilon > 0$ is arbitrary and $Y=-\Delta_f \ln H(x,t,y)$, we conclude that
\begin{equation}\label{eq19}
	\begin{aligned}
		t\left(\frac{|\nabla H|^2}{H^2}-\frac{\partial_t H}{H} \right) &\le \frac{n}{2}+C_{13}(L+\sqrt{K})\sqrt{(1+Kt)(1+t)}diam_M\\
		&\ \ \ + C_{13}(L+\sqrt{K})\sqrt{(1+Kt)(t+Kt+L^2t+A^2t+A^2Kt)},\ \forall t> 0.
	\end{aligned}
\end{equation}

From (\ref{eq19}), a short argument from \cite{ZQ} (see also \cite{ZFF3}) implies that the same bound actually holds if one replaces the $f$-heat kernel by any positive solution of the $f$-heat equation.

    This completes the proof of Theorem \ref{thm1}.                    \\ \qed
\section{ Li-Yau gradient bounds on closed manifolds under  $m$-Bakry-\'Emery Ricci curvature conditions  }
In this section, we prove Theorem \ref{thm2}. Since the proof of Theorem \ref{thm2} is similar, so we only present the key steps.

Before starting the proof of Theorem \ref{thm2}, let us present some results needed. First of all, Yi Li showed Li-Yau gradient estimate for manifolds with $m$-Bakry-\'Emery Ricci curvature bounded below.

\begin{lemma}\label{lem5} (\cite{Li})(Li-Yau gradient estimate)
	Let $(M^n,g)$ be a complete Riemannian manifold  with $Ric^{m,n}_f \ge -K$	for $K\ge0$. Then any positive solution $u$ of the f-heat equation $(\Delta_f -\partial_t)u=0$ on $M\times(0,T]$ satisfies 
	\begin{equation}
		\frac{|\nabla u|^2}{u^2}-\alpha\frac{u_t}{u}\le \frac{m\alpha^2K}{\alpha-1}+\frac{m\alpha^2}{2t}\nonumber
	\end{equation}
	for any $\alpha>1$.
\end{lemma}
	Based on the above Lemma, we can easily get the Harnack inequality.	
\begin{corollary}\label{cor2}(Harnack inequality)
	Under the same hypotheses as Lemma \ref{lem5}, we have
	\begin{equation}
		u(x,t_1) \le u(y,t_2)\left(\frac{t_2}{t_1}\right)^{\frac{m\alpha}{2}}e^{\left(\frac{\alpha d^2(x,y)}{4(t_2-t_1)}+\frac{m\alpha K}{\alpha-1}(t_2-t_1)\right)}\nonumber
	\end{equation}
	for all $x,y\in M$ and $0<t_1<t_2 \le T$.
\end{corollary}
In \cite{WWW}, Guofang Wei and Will Wylie proved the volume comparision for manifolds with $m$-Bakry-\'Emery Ricci curvature bounded below (see also \cite{DW}).
\begin{lemma} \label{lem6} (\cite{WWW})(volume comparison)
	Let $(M^n,g,e^{-f}dv)$ be a complete Riemannian manifold with $Ric^{m,n}_f \ge -(m-1)H$ for $H\ge0$. Then $\frac{V_f(B_p(R))}{V^m_H(B(R))}$ is nonincreasing in $R$.
\end{lemma}
In \cite{ALR}, Nelia Charalambous, Zhiqin Lu and Julie Rowlett got the upper and lower bounds for the $f$-heat kernel.
\begin{lemma} \label{lem7} (\cite{ALR})
	Let $(M^n,g,e^{-f}dv)$ be a complete Riemannian manifold with $Ric^{m,n}_f \ge -K$ on $B_{o}(4R+4) \subset M$ for $K \ge 0$ where $o$ is a point on $M$ and $R>0$. Then for any $x,y\in B_{o}(\frac{R}{4})$, $ 0<t<\frac{R^2}{4}$, and $\delta_1 \in (0,1)$
	\begin{equation}
		H(x,t,y)\ge \widetilde{c_6}(\delta_1,m)V^{-\frac{1}{2}}_f(B_x(\sqrt{t}))V^{-\frac{1}{2}}_f(B_y(\sqrt{t})) e^{\left(-\widetilde{c_7}(\delta_1,m)\frac{d^2(x,y)}{t}-\widetilde{c_8}(m)Kt\right)}\nonumber
	\end{equation}
	and 
	\begin{equation}
		H(x,t,y)\le \widetilde{c_3}(\delta_1,m)V^{-\frac{1}{2}}_f(B_x(\sqrt{t}))V^{-\frac{1}{2}}_f(B_y(\sqrt{t})) e^{\left(-\lambda_{1,f}(M)t-\frac{d^2(x,y)}{\widetilde{c_4}(\delta_1,m)t}+\widetilde{c_5}(m)\sqrt{Kt}\right)}\nonumber
	\end{equation}
	for some positive constants $\widetilde{c_3}(\delta_1,m),\ \widetilde{c_4}(\delta_1,m)$, \ $\widetilde{c_5}(m)$,\ $\widetilde{c_6}(\delta_1,m), \widetilde{c_7}(\delta_1,m)$ and $\widetilde{c_8}(m)$, where $\lambda_{1,f}(M)$ is the infimum of the weighted Rayleigh quotient on $M$.
	
	Whenever $Ric^{m,n}_f \ge -K$ on $M$ with $K\ge0$, then the same bound also holds for all $x,y\in M$ and $t>0$.
\end{lemma}
Now we are ready to prove Theorem \ref{thm2}.

\noindent{\it Proof of Theorem \ref{thm2}.}  
Let $u=u(x,t)$ be a positive solution of the $f$-heat equation on $M\times(0,+\infty)$. 
Write 
\begin{equation}
	Y=Y(x,t)=\left(\frac{|\nabla u|^2}{u^2}-\frac{\partial_t u}{u}\right)=-\Delta_f \ln u . \nonumber
\end{equation}
The Bochner formula \cite{Li} for  $Ric^{m,n}_f$ is 
\begin{equation} \label{eee3.1}
	\begin{aligned}
		\frac{1}{2}\Delta_f|\nabla u|^2&=|Hess\ u|^2+\left<\nabla \Delta_f u,\nabla u\right>+Ric^{m,n}_f(\nabla u,\nabla u)+\frac{1}{m-n}|\left<\nabla f,\nabla u\right>|^2\\
		&\ge\frac{(\Delta_f u)^2}{m}+\left<\nabla \Delta_f u,\nabla u\right>+Ric^{m,n}_f(\nabla u,\nabla u).
	\end{aligned}
\end{equation}

Using \eqref{eee3.1} and the $m$-Bakry-\'Emery Ricci curvature condition, we have 
\begin{equation} \label{eq20}
	\begin{aligned}
		(\Delta_f-\partial_t)Y+2\left<\nabla Y,\nabla \ln u\right>&=2|Hess \ln u|^2+2Ric_f^{m,n}(\nabla \ln u,\nabla \ln u)+\frac{2}{m-n} |\left<\nabla f,\nabla \ln u\right>|^2\\
		&\ge \frac{2(\Delta_f \ln u )^2}{m}-2K|\nabla \ln u |^2\\
		&=\frac{2}{m}Y^2-2K|\nabla \ln u |^2.
	\end{aligned}
\end{equation}
Let $Y^+(x,t)=\max\{Y(x,t),0\}$. Then inequality \eqref{eq20} implies that $Y^+$ is a subsolution of the inequality in the weak sense: on $M \times(0,+\infty)$,
\begin{equation} \label{eq21}
	(\Delta_f-\partial_t)Y^+ +2\left<\nabla Y^+,\nabla \ln u \right> \ge \frac{2}{m}(Y^+)^2-2K|\nabla \ln u|^2.
\end{equation}
For a positive integer $j$ and a small positive number $\epsilon$, multiplying both sides of \eqref{eq21} by $((t-\epsilon)^+)^{2j+2}(Y^+)^{2j}$ and integrating on $M\times (0,T]$ for any $T>0$ give that
\begin{equation}\label{eq22}
	\begin{aligned}
		&\frac{2}{m}\int_{0}^{T}\int_M ((t-\epsilon)^+Y^+)^{2j+2} d\mu dt \le \int_{0}^{T}\int_M ((t-\epsilon)^+)^{2j+2}(Y^+)^{2j}(\Delta_f -\partial_t)Y^+ d\mu dt\\
		&+ 2\int_{0}^{T}\int_M ((t-\epsilon)^+)^{2j+2}(Y^+)^{2j} \left<\nabla Y^+,\nabla \ln u \right>d\mu dt\\
		&+ 2K\int_{0}^{T}\int_M |\nabla \ln u|^2((t-\epsilon)^+)^{2j+2}(Y^+)^{2j}d\mu dt=T_1+T_2+T_3.
	\end{aligned}
\end{equation}
Using integration by parts as (\ref{eq4}) and (\ref{eq5}) for $T_1$ and $T_2$, we have the similar inequality

\begin{equation} \label{eq23}
	\begin{aligned}
		&(\frac{2}{m}-\frac{2}{2j+1})\int_{0}^{T}\int_M ((t-\epsilon)^+Y^+)^{2j+2}d\mu dt\\
		&\le \frac{2j+2}{2j+1}\int_{0}^{T}\int_M ((t-\epsilon)^+Y^+)^{2j+1} d\mu dt +\underbrace{2K\int_{0}^{T}\int_M |\nabla \ln u|^2((t-\epsilon)^+)^{2j+2}(Y^+)^{2j}d\mu dt}_{T_3}.
	\end{aligned} 
\end{equation} 

Using the upper and lower bounds (Lemma \ref{lem7}) for the $f$-heat kernel $H(x,t,y)$ and Harnack inequality (Corollary \ref{cor2}), Hamilton type estimate (Lemma \ref{lem1}) and the volume comparison theorem (Lemma \ref{lem6}) for manifolds with the $m$-Bakry-\'Emery Ricci curvature bounded below, we can find the bound for $t|\nabla_x \ln H(x,t,y)|^2$.

For a time $t_0>0$, we consider the $f$-heat kernel $H(x,t+t_0,y)$ with $t\in[0,t_0]$ at a fixed point $y\in M$.

According to Lemma \ref{lem7}, we choose $\delta =\frac{1}{2}$ and there exist some positive constants $\widetilde{C_1}$, $\widetilde{C_2}$, $\widetilde{C_3}$, $\widetilde{C_4}$, $\widetilde{C_5}$ and $\widetilde{C_6}$ depending only on $m$, such that
\begin{equation}\label{eq24}
	\frac{\widetilde{C_1}e^{[-\frac{d^2(x,y)}{\widetilde{C_2}(t+t_0)}+\widetilde{C_3}\sqrt{K(t+t_0)}]}}{V_f^{\frac{1}{2}}(B_x(\sqrt{t+t_0}))V_f^{\frac{1}{2}}(B_y(\sqrt{t+t_0}))} \ge H(x,t+t_0,y) \ge \frac{\widetilde{C_4}e^{[-\widetilde{C_5}\frac{d^2(x,y)}{(t+t_0)}-\widetilde{C_6}K(t+t_0)]}}{V_f^{\frac{1}{2}}(B_x(\sqrt{t+t_0}))V_f^{\frac{1}{2}}(B_y(\sqrt{t+t_0}))}.
\end{equation}
The upper bound implies 
\begin{equation}
	B:=\sup_{ M\times(0,t_0)}H(x,t+t_0,y)\le \frac{\widetilde{C_7}e^{\widetilde{C_8}Kt_0}}{\inf_{z\in M}V_f(B_z(\sqrt{t_0}))},\nonumber
\end{equation}  
which yields, by the lower bound of $H(x,t+t_0,y)$, that
\begin{equation}\label{eqq25}
	\frac{B}{H(x,t+t_0,y)}\le \widetilde{C_9}e^{[\widetilde{C_{10}}Kt_0+\widetilde{C_5}\frac{d^2(x,y)}{t_0}]} 
	\frac{\sup_{w\in M}V_f(B_w(\sqrt{2t_0}))}{\inf_{z\in M}V_f(B_z(\sqrt{t_0}))}.
\end{equation} 
We notice the infimum and supremum of the geodesic ball's volumes can be achieved by some points on $M$, say point $p$ and point $q$, i.e.,
\begin{equation}\label{eq26}
	\frac{\sup_{w\in M}V_f(B_w(\sqrt{2t_0}))}{\inf_{x\in M}V_f(B_z(\sqrt{t_0}))}=\frac{V_f(B_q(\sqrt{2t_0}))}{V_f(B_p(\sqrt{t_0}))}.
\end{equation}
By lemma \ref{lem6}, we get
\begin{equation}\label{eq27}
	\begin{aligned}
	\frac{V_f(B_p(\sqrt{2t_0}+d(p,q)))}{V_f(B_p(\sqrt{2t_0}))}
	&\le \frac {\int_{0}^{d(p,q)+\sqrt{2t_0}}\left(sinh(\sqrt{\frac{K}{m-1}}r)\right)^{m-1}dr}{\int_{0}^{\sqrt{2t_0}}\left(sinh(\sqrt{\frac{K}{m-1}}r)\right)^{m-1}dr}\\ &\le \left(\frac{d(p,q)}{\sqrt{2t_0}}+1\right)^me^{\sqrt{K(m-1)}d(p,q)}\\
	&=e^{m\ln \left(\frac{d(p,q)}{\sqrt{2t_0}}+1\right) +\sqrt{K(m-1)}d(p,q)}\\
	&\le e^{m  \frac{d(p,q)}{\sqrt{2t_0}}+ \sqrt{K(m-1)}d(p,q)},  
\end{aligned}
\end{equation}
and 
\begin{equation}
	\frac{V_f(B_p(\sqrt{2t_0}))}{V_f(B_p(\sqrt{t_0}))} \le \widetilde{C_{12}}e^{\widetilde{C_{11}}\sqrt{Kt_0}    } . \nonumber
\end{equation}
This and \eqref{eqq25} imply that
\begin{equation}
	\ln \frac{B}{H(x,t+t_0,y)} \le \ln \widetilde{C_9} +
	\widetilde{C_{10}}Kt_0+\widetilde{C_5}\frac{diam_M^2}{t_0}+\ln \widetilde{C_{12}}+
	\widetilde{C_{11}}\sqrt{Kt_0}+\frac{mdiam_M}{\sqrt{2t_0}}+\sqrt{K(m-1)}diam_M.\nonumber
\end{equation}
Since the $m$-Bakry-\'Emery Ricci curvature bounded below, we can deduce Bakry-\'Emery Ricci curvature bounded below, so we can still use Lemma \ref{lem1} and that
\begin{equation}
	t|\nabla_x \ln H(x,t+t_0,y)|^2 \le\widetilde{C_{13}} (1+2Kt_0)\left(1+Kt_0+
	\frac{diam^2_M}{t_0}\right), t\in (0,t_0].\nonumber
\end{equation}
Then we conclude
\begin{equation} \label{eq26}
	\begin{aligned}
		t|\nabla_x \ln H(x,t,y)|^2 \le\widetilde{C_{14}} (1+Kt)\left(1+Kt+
		\frac{diam^2_M}{t}\right),\ \forall t>0.
	\end{aligned}
\end{equation}

The bound is adequate for us when then time is short, say $t \le 4$. For any time $t\ge 4$, since $\int_M H(x,t+1,y)d\mu(x)=1$, there is a point $x_1\in M$ such that  $H(x_1,t+1,y)=\frac{1}{V_f(M)}$. According to Corollary \ref{cor2} with $\alpha=2,\ t_1=t,\ t_2=t+1$, we have     
\begin{equation} 
	H(x,t,y) \le H(x_1,t+1,y)\left(\frac{t+1}{t}\right)^m e^{\left(2mK+\frac{d^2(x,x_1)}{2} \right)}. \nonumber
\end{equation}  
Since $t\ge4$, this implies
\begin{equation} \label{eq29}
       	H(x,t,y) \le \left(\frac{5}{4}\right)^m \frac{1}{V_f(M)} e^{\left(2mK+\frac{diam^2_M}{2}\right)}:=B.
\end{equation}                
Similarly, there is a point $x_2$ such that $H(x_2,t-1,y)=\frac{1}{V_f(M)}$ and that
\begin{equation} 
	H(x_2,t-1,y) \le H(x,t,y)\left(\frac{t}{t-1}\right)^m e^{\left(2mK+\frac{d^2(x,x_2)}{2}\right)}, \nonumber
\end{equation}  
which implies
\begin{equation} \label{eq30}
	H(x,t,y) \ge \left(\frac{3}{4}\right)^m \frac{1}{V_f(M)} e^{\left(-2mK-\frac{diam^2_M}{2}\right)}, \ t\ge4.
\end{equation}  
Using (\ref{eq29}) and (\ref{eq30}), we find, for $t_0 \ge4$, that
\begin{equation}
	\ln \frac{B}{H(x,t+t_0,y)} \le m\ln2+4mK+diam^2_M,\ t\in(0,t_0].\nonumber
\end{equation}
This and Lemma \ref{lem1} yield
\begin{equation}
	t|\nabla_x \ln H(x,t+t_0,y)|^2\le(1+2Kt_0)(m\ln2+4mK+diam^2_M),\ t\in(0,t_0].\nonumber
\end{equation}
Therefore
\begin{equation} \label{eq31}
	t|\nabla_x \ln H(x,t,y)|^2 \le2(1+Kt)(m\ln2+4mK+diam^2_M),\ t\ge4.
\end{equation}
Next, plugging (\ref{eq26}) for $t<4$ and (\ref{eq31}) for $t\ge4$ into the term $T_3$ in (\ref{eq23}) with $u=H(x,t,y)$, we obtain 
\begin{equation}
T_3 \le K(1+KT)(\widetilde{C_{15}}(T+KT+diam^2_M)+Tdiam^2_M)\int_{0}^{T}\int_M ((t-\epsilon)^+Y^+)^{2j}d\mu dt.\nonumber
\end{equation}
These are the key steps to prove Theorem \ref{thm2} and the rest of the proof is similar to Theorem \ref{thm1}. \\ \qed

\section*{Acknowledgements}

Research is partially supported by NSFC Grant No. 11971168, Shanghai Science and Technology Innovation Program Basic Research Project STCSM 20JC1412900.


\begin{thebibliography}{99}

\bibitem{BE} Bakry, Dominique; Émery, Michel, \emph{Diffusions hypercontractives}. [Hypercontractive diffusions] Séminaire de probabilités, XIX, 1983/84, 177–206. 
\bibitem{BQ}Bakry, Dominique; Qian, Zhongmin, \emph{Harnack inequalities on a manifold with positive or negative Ricci curvature}. Rev. Mat. Iberoamericana 15 (1999), no. 1, 143–179.
\bibitem{Bam} Bamler, Richard, \emph{Convergence of Ricci flows with bounded scalar curvature.}  Ann. of Math. (2) 188 (2018), no. 3, 753-831.
\bibitem{ALR} Charalambous, Nelia; Lu, Zhiqin;  Rowlett, Julie, \emph{Eigenvalue estimates on Bakry-\'Emery manifolds}. Elliptic and parabolic equations, 45-61, Springer Proc. Math. Stat., 119, Springer, Cham, 2015.
\bibitem{CW} Chen, Xiuxiong; Wang, Bing,  \emph{Space of Ricci flows (II)-Part B: Weak compactness of the flows.} J. Differential Geom. 116 (2020), no. 1, 1-123.
\bibitem{BLN} Chow, Bennett; Lu, Peng; Ni, Lei, \emph{Hamilton’s Ricci Flow}, Graduate Studies in Mathematics, vol. 77, American Mathematical Society/Science Press, Providence, RI/New York, 2006.
\bibitem{DW} Dai, Xianzhe; Wei, Guofang, \emph{Comparison Geometry for Ricci Curvature}.

http://web.math.ucsb.edu/dai/preprints.html
\bibitem{GM} Garofalo, Nicola; Mondino, Andrea, \emph{Li–Yau and Harnack type inequalities in RCD∗(K; N) metric measure spaces}, Nonlinear Anal. 95 (2014) 721–734.
\bibitem{HA} Hamilton, Richard S., \emph{A matrix Harnack estimate for the heat equation}.Comm. Anal. Geom. 1 (1993), no. 1, 113–126.
\bibitem{LX} Li, Junfang; Xu, Xiangjin, \emph{Differential Harnack inequalities on Riemannian manifolds I: linear heat equation}, Adv. Math. 226 (5) (2011) 4456–4491.
\bibitem{LY} Li, Peter; Yau, Shing-Tung, \emph{On the parabolic kernel of the Schrödinger operator}.
Acta Math. 156 (1986), no. 3-4, 153–201.
\bibitem{Li}Li, Yi, \emph{Li-Yau-Hamilton estimates and Bakry-Emery-Ricci curvature}.Nonlinear Anal. 113 (2015), 1–32.
\bibitem{Lo} Lott, John, \emph{ Some geometric properties of the Bakry–Émery–Ricci tensor}. Comment. Math. Helv. 78,865–883 (2003) 
\bibitem{Per} Perelman, Grisha, \emph{The entropy formula for the Ricci flow and its geometric applications,} arXiv:math/0211159.
\bibitem{Qian Bin}Qian, Bin, \emph{Remarks on differential Harnack inequalities}. J. Math. Anal. Appl. 409 (2014), no. 1, 556–566.
\bibitem{QZZ} Qian, Zhongmin; Zhang, Huichun; Zhu, Xiping, \emph{Sharp spectral gap and Li–Yau’s estimate on Alexandrov spaces}, Math. Z. 273 (3–4) (2013) 1175–1195.
\bibitem{TZ} Tian, Gang; Zhang, Zhenlei, \emph{Regularity of K\"ahler-Ricci flows on Fano manifolds.} Acta Math. 216 (2016), no. 1, 127-176.
\bibitem{WFY}Wang, Fengyu, \emph{Gradient and Harnack inequalities on noncompact manifolds with boundary}, Pacific J. Math. 245 (1) (2010) 185–200.
\bibitem{WWW}Wei, Guofang; Wylie, Will, \emph{ Comparison geometry for the Bakry–Émery Ricci tensor}. J. Differ. Geom. 83,377–405 (2009)
\bibitem{WJY3} Wu, Jiayong; Wu, Peng, \emph{$L^p$-Liouville theorems on complete smooth metric measure spaces}. Bull. Sci. Math. 138,510–539 (2014)
\bibitem{WJY2} Wu, Jiayong; Wu, Peng, \emph{Heat kernel on smooth metric measure spaces with nonnegative curvature}. Math.Ann. 362, 717–742 (2015)
\bibitem{WWu}Wu, Jiayong; Wu, Peng, \emph{Heat kernel on smooth metric measure spaces and applications}. Math. Ann. 365 (2016), no. 1-2, 309–344.
\bibitem{Yau}Yau, Shing-Tung, \emph{Harnack inequality for non-self-adjoint evolution equations}.
Math. Res. Lett. 2 (1995), no. 4, 387–399.
\bibitem{ZFF2} Yu, Chengjie; Zhao, Feifei, \emph{Sharp Li-Yau-type gradient estimates on hyperbolic spaces}. J. Geom. Anal. 30 (2020), no. 1, 54–68.
\bibitem{ZFF} Yu, Chengjie; Zhao, Feifei, \emph{Recurrence relations for heat kernels on spheres}. J. Math. Anal. Appl. 507 (2022), no. 2, Paper No. 125790, 14 pp.
\bibitem{ZFF3} Yu, Chengjie; Zhao, Feifei, \emph{Li-Yau multiplier set and optimal Li-Yau gradient estimate on hyperbolic spaces}. Potential Anal. 56 (2022), no. 2, 191–211.
\bibitem{ZZH} Zhang, Huichun; Zhu, Xiping, \emph{Local Li–Yau’s estimates on RCD∗(K, N) metric measure spaces}, Calc. Var. Partial Differential Equations (2016) 55–93.
\bibitem{ZQ}Zhang, Qi S., \emph{A Sharp Li-Yau gradient bound on Compact Manifolds }, arXiv:2110.08933.









\bibitem{ZZ} Zhang, Qi S.; Zhu, Meng, \emph{Li-Yau gradient bound for collapsing manifolds under integral curvature condition}. Proc. Amer. Math. Soc. 145 (2017), no. 7, 3117–3126.
\bibitem{ZZ2} Zhang, Qi S.; Zhu, Meng, \emph{Li-Yau gradient bounds on compact manifolds under nearly optimal curvature conditions}. J. Funct. Anal. 275 (2018), no. 2, 478–515.







\end{thebibliography}
\end{document}